\documentclass[12pt]{article}
\usepackage{latexsym}

\usepackage[tbtags]{amsmath}
\usepackage{epsfig,amstext,amssymb,amsthm,latexsym}
\usepackage{amssymb,color}

\definecolor{c20}{rgb}{0.,0.7,0.}
\definecolor{c30}{rgb}{0.,0.,1.}
\definecolor{c40}{rgb}{1,0.1,0.7}
\definecolor{c50}{rgb}{1,0,0}

\date{}
\newtheorem{theorem}{Theorem}[section]
\newtheorem{lemma}{Lemma}[section]

\newtheorem{remark}{Remark}[section]
\newtheorem{example}{Example}[section]

\newtheorem{property}{Property}[section]

\makeatletter 
\@addtoreset{equation}{section}
\makeatother 

\begin{document}

\title{On the Kendall Correlation Coefficient}

\author{ A. Stepanov, \thanks{\noindent Immanuel Kant Baltic Federal University, A.Nevskogo 14, Kaliningrad, 236041 Russia; 
}
~\small{\it Immanuel Kant Baltic Federal University} }
\maketitle
\begin{abstract} In the present paper, we  first discuss the Kendall rank correlation coefficient $\tau_n$. In continuous case, we define $\tau_n$ in terms of the concomitants of order statistics, find the expected value of $\tau_n$ and show that the later is free of $n$. We also prove that in  continuous case the Kendall  correlation coefficient converges in probability to its expected value $\tau=E\tau_n$. We then propose to consider $\tau$ as a new theoretical correlation coefficient which can be an alternative to the classical  Pearson product-moment correlation coefficient. At the end of this work we analyze illustrative examples.
\end{abstract}
\noindent {\it Keywords and Phrases}:  bivariate distributions;
 concomitants of order statistics;    Pearson product-moment correlation coefficient; sample correlation coefficient; Kendall rank correlation coefficient.

\section{ Introduction}
Let  $(X,Y), (X_1,Y_1), (X_2,Y_2),\ldots,(X_n,Y_n)$
be independent  and identically distributed random vectors with
bivariate distribution $F(x,y)=P(X<x,Y<y)$ and  corresponding
marginal distributions $H(x)=P(X<x)$ and $G(y)=P(Y<y)$. The purpose of  slight modifications in the definitions of $F, H$ and $G$ is that in Section~4 we define a new correlation coefficient and we wish to have the same form of this correlation coefficient for  different types of distributions. If  $F$ is an absolutely continuous distribution then the corresponding densities will be denoted as $f(x,y), h(x)$ and $g(y)$, respectively. Let  $X_{1,n}\leq X_{2,n}\leq\ldots\leq X_{n,n}$ be the order statistics obtained from the sample  $X_1,X_2,\ldots,X_n$. For these order statistics, let us define their concomitants $Y_{[1,n]},Y_{[2,n]},\ldots,Y_{[n,n]}$. Let $X_i=X_{j,n}$. Then $Y_{[j,n]}=Y_i$ is the concomitant of the order statistic $X_{j,n}$.  Concomitants of order statistics were proposed by David (1973) and Bhattacharya (1974). Concomitants of order statistics were further discussed  in  David and Galambos (1974),  Bhattacharya (1984), Egorov and Nevzorov (1984),  David (1994),     Goel and Hall (1994), Chu {\it et al.} (1999), David and Nagaraja (2003), Bairamov and Stepanov (2010) and others.

It is known that the rate of dependence between random variables $X$ and $Y$ can be measured in terms of  the Pearson product-moment correlation coefficient
$$
\rho={E(X-EX)(Y-EY) \over \sigma_X\sigma_Y}.
$$
The sample correlation coefficient 
$$\rho_n =\frac{\sum\limits_{i=1}^n (X_i-\bar{X})(Y_i-\bar{Y})} {\sqrt{\sum\limits_{i=1}^n (X_i-\bar{X})^2 \sum\limits_{i=1}^n (Y_i-\bar{Y})^2}}
$$
is a good approximation for $\rho $ for large values of $n$ since $\rho_n\stackrel{p}{\rightarrow}\rho$; see, for example, Fisher (1921).

Let   $Y_{[1,n]},\ldots, Y_{[n,n]}$ be the  concomitants of order statistics $X_{1,n}\leq \ldots\leq  X_{n,n}$. Let
$$
r_{j,n}=\sum_{i=1}^{j-1}I_{i,j,n}=\sum_{i=1}^{j-1}I_{\{Y_{[i,n]}\leq Y_{[j,n]}\}}\quad (1\leq i<j\leq n)
$$
be the rank of the concomitant $Y_{[j,n]}$ among the concomitants $Y_{[1,n]},\ldots, Y_{[j,n]}$. The following value
$$
\tau _n=\frac{4\sum_{j=2}^nr_{j,n}}{n(n-1)}-1=\frac{4}{n(n-1)}\sum_{j=2}^n\sum_{i=1}^{j-1}I_{i,j,n}-1
$$
is known as Kendall's rank correlation coefficient (see, for example, Kendall (1970)). This rank correlation coefficient $\tau_n$ can be used along with $\rho_n$. Sometimes  samples $X_i$ and $Y_i$ are not known but the ranks $r_{j,n}$  are known. In this case one can use $\tau_n$ instead of $\rho_n$. It should be noted that there is another rank correlation coefficient -- Spearman's rank correlation coefficient, which we do not discuss in the present work.

Basic properties of $\tau_n$ are as follows. If the agreement/disagreement between the sequence $X_{j,n}$ and the rankings of $Y_{[j,n]}$ is perfect,  the coefficient  value $\tau_n$ is near 1/-1. Further, if  variables $X$ and $Y$ are independent, then  $Y_{[j,n]}\ j=1,\ldots,n$ are  independent and identically distributed. Obviously, $P(Y_{[i,n]}\leq Y_{[j,n]})=1/2\ \forall i, j, i\not= j$. Then $\tau _n\stackrel{a.s.}{\rightarrow}0$. Different bounds for $\tau_n$ (basically in the normal case) were obtained in Daniels (1950), Durbin and Stuart (1951), Kendall (1970) and Xu {\it et al.} (2009); see also references in these works.

Further in the paper, we discuss moment and asymptotic properties of $\tau_n$. In Section~2, we find  $E\tau_n$ in  general continuous case and show that it is free of $n$. In Section~3, we prove that in  general continuous case the Kendall  correlation coefficient converges in probability to its expected value. The last observation motivates us to introduce in Section~4 a new correlation coefficient: $\tau=E\tau_n$. This correlation coefficient $\tau$ is  a possible alternative to the classical  Pearson product-moment correlation coefficient $\rho$. It should be noted that  no moment assumption is needed for the existence of $\tau$. We illustrate our theoretical results in Section~5 by examples and simulation results. 

\section{The expected value of $\tau_n$}
We assume in this section that $F$ is a continuous distribution.
\begin{lemma}\label{lemma2.1} For $n\geq 2$, we have
\begin{eqnarray}\label{2.1}
E\tau _n&=&4\int_{x\leq u, y\leq v}F(dx,dy)F(du,dv)-1\\
\label{2.2}&=&4\int_{\mathbb{R}^2}(1-H(x)-G(y)+F(x,y))F(dx,dy)-1\\
\label{2.3}&=&4\int_{\mathbb{R}^2}F(x,y)F(dx,dy)-1=4EF(X,Y)-1.
\end{eqnarray}
\end{lemma}
\begin{gproof}{} It is easily seen that (\ref{2.1}) implies both (\ref{2.2}) and (\ref{2.3}). We present the proof (\ref{2.1}) for the case when $F$ is an absolutely continuous distribution.  Let  
$$
Z_{[1,n]}=(X_{1,n},Y_{[1,n]}),\ldots,Z_{[n,n]}=(X_{n,n},Y_{[n,n]}).
$$
One can find   that
\begin{equation}\label{2.4}
 f_{Z_{[1,n]},\ldots,Z_{[n,n]}}(x_1,y_1,\ldots,x_n,y_n)=n!f(x_1,y_1)\ldots f(x_n,y_n)\quad (x_1\leq \ldots\leq x_n, y_i\in\mathbb{R} )
\end{equation}
and
\begin{equation}\label{2.5}
 f_{Y_{[1,n]},Y_{[2,n]},\ldots,Y_{[n,n]}}(y_1,y_2,\ldots,y_n)=n!\int_{x_1\leq \ldots\leq x_n} f(x_1,y_1)\ldots f(x_n,y_n)dx_1\ldots dx_n.
\end{equation}
By integration in (\ref{2.5}), one can find  that for $y,v\in \mathbb{R}$
\begin{eqnarray*}
f_{Y_{[i,n]},Y_{[j,n]}}(y,v)&=&\frac{n!}{(i-1)!(j-i-1)!(n-j)!}\\
&\times &\int_{x\leq u}H^{i-1}(x) f(x,y)(H(u)-H(x))^{j-i-1}f(u,v)(1-H(u))^{n-j}dxdu.
\end{eqnarray*}
Then
\begin{eqnarray*}
Er_{j,n}&=&\sum_{i=1}^{j-1}P(Y_{[i,n]}\leq Y_{[j,n]})=\frac{n!}{(j-2)!(n-j)!}\\
&\times &\int_{x\leq u, y\leq v}f(x,y)f(u,v)H^{j-2}(u) (1-H(u))^{n-j}dxdudydv.
\end{eqnarray*}
We finish at the following identity
\begin{equation}\label{2.6}
E\tau _n=\frac{4\sum_{j=2}^nEr_{j,n}}{n(n-1)}-1=4\int_{x\leq u, y\leq v}f(x,y)f(u,v)dxdydudv-1.
\end{equation}
The result readily follows. When $F$ is continuous the result can be proved in the same manner.
\end{gproof}
\begin{remark}\label{remark2.1}
Since $E\tau_n$ is free of $n$, let $\tau=E\tau_n$.
\end{remark}

\section{Asymptotic properties of $\tau_n$}
We assume in this section that $F$ is a continuous distribution.
\begin{theorem}\label{theorem3.1} The following asymptotic property holds true
$$
\tau_n\stackrel{p}{\rightarrow}\tau.
$$
\end{theorem}
\begin{gproof}{} We present the proof of Theorem~\ref{theorem3.1} for the case when $F$ is an absolutely continuous distribution. It follows from Chebyshev's inequality that for any $\varepsilon >0$
\begin{eqnarray}
\nonumber&&P(\mid \tau_n-\tau\mid >\varepsilon )\leq \frac{Var \tau_n}{\varepsilon ^2}\\
&&=\frac{16}{\varepsilon ^2n^2(n-1)^2}\left[\sum_{j=2}^n\sum_{i=1}^{j-1}\sum_{k=2}^n\sum_{l=1}^{k-1}E(I_{i,j,n}I_{l,k,n})-\left(E\sum_{j=2}^n\sum_{i=1}^{j-1}I_{i,j,n}\right)^2\right].\label{3.1}
\end{eqnarray}
By (\ref{2.6}), we have 
$$
\left(E\sum_{j=2}^n\sum_{i=1}^{j-1}I_{i,j,n}\right)^2=\left(n(n-1)\int_{x\leq u, y\leq v}f(x,y)f(u,v)dxdydudv\right)^2.
$$
Let 
$$
\sum=\sum_{j=2}^n\sum_{i=1}^{j-1}\sum_{k=2}^n\sum_{l=1}^{k-1}E(I_{i,j,n}I_{l,k,n})=\sum_{j=2}^n\sum_{i=1}^{j-1}\sum_{k=2}^n\sum_{l=1}^{k-1}P(Y_{[i,n]}\leq  Y_{[j,n]},Y_{[l,n]}\leq Y_{[k,n]}).
$$
One can write $\Sigma=\Sigma_1+\Sigma_2$, where
\begin{eqnarray*}
\Sigma_1&=&\sum_{1\leq i<j<l<k\leq n}+\sum_{1\leq i<l<j<k\leq n}+\sum_{1\leq l<i<j<k\leq n}\\
\nonumber &&\sum_{1\leq l<k<i<j\leq n}+\sum_{1\leq l<i<k<j\leq n}+\sum_{1\leq i<l<k<j\leq n},
\end{eqnarray*}
\begin{eqnarray*}
\Sigma_2&=&2\sum_{1\leq i<j=l<k\leq n}+2\sum_{1\leq i=l<j<k\leq n}+2\sum_{1\leq i<l<j=k\leq n}\\
 &&+4\sum_{1\leq i=l<j=k\leq n}+2\sum_{1\leq l<i<j=k\leq n}+2\sum_{1\leq l<k=i<j\leq n}+2\sum_{1\leq l=i<k<j\leq n}.
\end{eqnarray*}
Here 
$$
\sum_{1\leq i<j<l<k\leq n}=\sum_{k=4}^n\sum_{l=3}^{k-1}\sum_{j=2}^{l-1}\sum_{i=1}^{j-1}P(Y_{[i,n]}\leq Y_{[j,n]},Y_{[l,n]}\leq  Y_{[k,n]}),
$$ 
$$
\sum_{1\leq i<j=l<k\leq n}=\sum_{k=3}^n\sum_{j=2}^{k-1}\sum_{i=1}^{j-1}P(Y_{[i,n]}\leq  Y_{[j,n]},Y_{[j,n]}\leq  Y_{[k,n]})
$$
and the other terms (sums) in $\Sigma_1$ and $\Sigma_2$ are designated in the same fashion. We show in this section that the terms in $\Sigma_1$ behave like $O(n^4)$ when the terms in $\Sigma_2$ behave like $o(n^4)$. Let us start with $\Sigma_2$.  Let us, for example,  show  that $\frac{\sum_{1\leq i<j=l<k\leq n}}{n^4}\rightarrow 0$ as $n\rightarrow \infty $. It follows from (\ref{2.5}) that 
\begin{eqnarray*}
&&\sum_{k=3}^n\sum_{j=2}^{k-1}\sum_{i=1}^{j-1}P(Y_{[i,n]}\leq  Y_{[j,n]}\leq  Y_{[k,n]})=\sum_{k=3}^n\sum_{j=2}^{k-1}\sum_{i=1}^{j-1}\frac{n!}{(i-1)!(j-i-1)!(k-j-1)!(n-k)!} \\
&& \times \int_{x_i\leq x_j\leq x_k, y_i\leq y_j\leq y_k}H^{i-1}(x_i)f(x_i,y_i)(H(x_j)-H(x_i))^{j-i-1}f(x_j,y_j)(H(x_k)-H(x_j))^{k-j-1}\\
&&\times f(x_k,y_k)(1-H(x_k))^{n-k}dx_idy_idx_jdy_jdx_kdy_k\\
&&=n(n-1)(n-2)\int_{x_i\leq x_j\leq x_k, y_i\leq y_j\leq y_k}f(x_i,y_i)f(x_j,y_j)f(x_k,y_k)dx_idy_idx_jdy_jdx_kdy_k\\
&&\leq n(n-1)(n-2)=o(n^4).
\end{eqnarray*}
One can show that the other terms of $\Sigma_2$ behave like $o(n^4)$ too. Let us now estimate  the  terms in $\Sigma_1$. Let us take, for example,  $\sum_{1\leq i<j<l<k\leq n}$. It follows from (\ref{2.5}) that 
\begin{eqnarray*}
&& \sum_{1\leq i<j<l<k\leq n}=\sum_{k=4}^n\sum_{l=3}^{k-1}\sum_{j=2}^{l-1}\sum_{i=1}^{j-1}P(Y_{[i,n]}\leq Y_{[j,n]}, Y_{[l,n]}\leq Y_{[k,n]})\\
&& =\sum_{k=4}^n\sum_{l=3}^{k-1}\sum_{j=2}^{l-1}\sum_{i=1}^{j-1}\frac{n!}{(i-1)!(j-i-1)!(l-j-1)!(k-l-1)!(n-k)!}\\
&&\times \int_{x_i\leq x_j\leq x_l\leq x_k, y_i\leq y_j, y_l\leq y_k}H^{i-1}(x_i)f(x_i,y_i)(H(x_j)-H(x_i))^{j-i-1}f(x_j,y_j)(H(x_l)-H(x_j))^{l-j-1}\\
&&\times f(x_l,y_l)(H(x_k)-H(x_l))^{k-l-1}f(x_k,y_k)(1-H(x_k))^{n-k}dx_idy_ix_jdy_jdx_ldy_ldx_kdy_k\\
&&=n(n-1)(n-2)(n-3)\int_{x_i\leq x_j\leq x_l\leq x_k, y_i\leq y_j, y_l\leq y_k}T,
\end{eqnarray*}
where 
$$
T=f(x_i,y_i)f(x_j,y_j) f(x_l,y_l)f(x_k,y_k)dx_idy_ix_jdy_jdx_ldy_ldx_kdy_k.
$$
In the same way one can estimate the other terms in $\Sigma_1$. Then
\begin{eqnarray*}
&&\sum=n(n-1)(n-2)(n-3) \\
&& \times \left(\int_{x_i\leq x_j\leq x_l\leq x_k, y_i\leq y_j, y_l\leq y_k}T+ \int_{x_i\leq x_l\leq x_j\leq x_k, y_i\leq y_j, y_l\leq y_k}T+\int_{x_l\leq x_i\leq x_j\leq x_k, y_i\leq y_j, y_l\leq y_k}T \right.\\
&&\left.+ \int_{x_l\leq x_k\leq x_i\leq x_j, y_i\leq y_j, y_l\leq y_k}T+ \int_{x_l\leq x_i\leq x_k\leq x_j, y_i\leq y_j, y_l\leq y_k}T+  \int_{x_l\leq x_k\leq x_i\leq x_j, y_i\leq y_j, y_l\leq y_k}T\right)+o(n^4)\\
&&=\int_{x_i\leq x_j, x_l\leq x_k, y_i\leq y_j, y_l\leq y_k}T+o(n^4).
\end{eqnarray*}
Obverve that
\begin{eqnarray*}
&&\int_{x_i\leq x_j, x_l\leq x_k,  y_i\leq y_j, y_l\leq y_k}T \\
&& =\int_{x_i\leq x_j,x_l\leq x_k, y_i\leq y_j, y_l\leq y_k}f(x_i,y_i)f(x_j,y_j) f(x_l,y_l)f(x_k,y_k)dx_idy_ix_jdy_jdx_ldy_ldx_kdy_k\\
&&=\left(\int_{x_i\leq x_j, y_i\leq y_j}f(x_i,y_i)f(x_j,y_j)dx_idy_ix_jdy_jdx\right)^2.
\end{eqnarray*}
It follows from (\ref{3.1}) that
$$
P(\mid \tau_n-\tau\mid >\varepsilon )\rightarrow 0\quad (n\rightarrow \infty).
$$
The result readily follows.  When $F$ is continuous the result can be proved in the same manner.
\end{gproof}

\section{New correlation coefficient} In  continuous case, we propose to consider the value 
$$
\tau=4EF(X,Y)-1=4\int_{x\leq u, y\leq v}F(dx,dy)F(du,dv)-1=4\int_{\mathbb{R}^2}F(x,y)F(dx,dy)-1
$$ 
as a new theoretical correlation coefficient, which measures  the rate of dependence between the random variables $X$ and $Y$. Basic properties of $\tau$ are as follows.

\begin{property}\label{property4.1}
It is obvious that $-1\leq \tau\leq 1$.
\end{property}
\begin{property}\label{property4.2} 
It is easily seen that if $X$ and $Y$ are independent, then $\tau=0$. 
\end{property}
\begin{property}\label{property4.3}
Let $Y=\varphi(X)$, where $\varphi$ is a  nondecreasing function. Then $\tau=1$. Let $Y=\psi (X)$, where $\psi $ is a  non increasing function. Then $\tau=-1$.
\end{property} 
\begin{gproof}{} 
We prove only the first statement. The second statement can be proved in the same manner. We again assume that $F$ is an absolutely continuous distribution. Let us consider the integral $\int_{\mathbb{R}^2}F_{X,Y}(x,y)f_{X,Y}(x,y)dxdy$, where the random variables $X$ and $Y$ are now attached to the corresponding bivariate distribution and density. It is obvious that 
$$
P(x\leq X< x+dx,y\leq Y< y+dy)\approx f_{X,Y}(x,y)dxdy.
$$
Let $Y=X$. Then for all small enough $dx$ and $dy$
$$
P(x\leq X< x+dx,y\leq X< y+dy)=\left\{
\begin{array}{cc}
h_X(x)\min\{dx,dy\}, & \mbox{if}\quad y=x\\
0, & \mbox{otherwise}.
\end{array}
\right.
$$
It follows from the last identity that
$$
\int_{\mathbb{R}^2, Y=X}F_{X,Y}(x,y)f_{X,Y}(x,y)dxdy=\int_{\mathbb{R}}H_{X}(x)h_{X}(x)dx=1/2.
$$
That way, if $Y=X$, then $\tau=4EF(X,Y)-1=1$. The result readily follows since $EF(X,\varphi (X))\geq EF(X,X)$.
\end{gproof}
Observe that the corresponding  Property~\ref{property4.3} for $\rho$ is as follows:  $\rho=\pm 1$ iff $Y=aX+b$.

In  discrete case, one can define  $\tau$ by the identity
$$
\tau=4EF(X,Y)-1=4\sum_j\sum_kf(j,k)F(j,k)-1,
$$
where $f(j,k)=P(X=j,Y=k)$ and $F(j,k)=P(X<j,Y<k)$. One can check here the validity of Properties~\ref{property4.1}-\ref{property4.3}. To check Property~\ref{property4.1} one can apply the inequality  $\sum_jh(j)H(j)\leq 1/2$, where $h(j)=P(X=j)$ and $H(j)=P(X<j)$.

In  general case, when $F$ has discrete and continuous components, one can define $\tau$ as follows: $\tau=4EF(X,Y)-1$.\\

The proposed correlation coefficient $\tau$ has some advantages and disadvantages in comparison with   the Pearson product-moment correlation coefficient $\rho$. 

Advantages.

(1)  We think  that  $\tau$ reflects the rate of dependence between  $X$ and $Y$ better than $\rho$ since $\tau$ is based on the whole information about a distribution. Observe that $\rho$ is based only on the first and second moments.

(2) No moment assumption is needed for the existence of $\tau$. 

(3) In  continuous case, we have  proved  that $\tau_n\stackrel{p}{\rightarrow}\tau$. It is also true that $\rho_n\stackrel{p}{\rightarrow}\rho$. That way, both $\tau$ and $\rho$ are approximated for large values of $n$ by  $\tau_n$ and $\rho_n$, respectfully. However, $E\tau_n=\tau$.

Disadvantages.

(1) The Pearson product-moment correlation coefficient
$\rho$ is simpler than $\tau$. To compute $\rho$ one should  only know the first and the second moments. To compute $\tau$  one should know the distribution $F$. In this respect, $\rho$ has an advantage over $\tau$, because in many real experiments (say,  in the time-series analysis) one  often works only with first and second moments. 

(2) The Pearson product-moment correlation coefficient
$\rho$ is visibly presented in many important statistical models and formulas such as linear regression models, bivariate normals densities and so on.

\section{Examples}
\begin{example}\label{example5.1}
Let
$$
F(x,y)=1-e^{-x}-\frac{1-e^{-x(y+1)^t}}{(y+1)^t}\quad (x>0, y>0, t>0)
$$
be a bivariate distribution with marginal distributions  $H(x)=1-e^{-x}\ (x>0)$ and
$G(y)=1-\frac{1}{(y+1)^t}\ (y>0)$.  It follows from (\ref{2.2}) that for any $t>0$
$$
\tau=E\tau_n=-0.5.
$$ 
We made a corresponding simulation experiment, i.e. $\tau_n$ was computed many times for "large" values of $n$ by simulation in Matlab. The code for Kendall's correlation coefficient is $$corrcoef(x,y,'type',Kendall).$$  The experiment gave us the same value $-0.5$. 

Observe that here $\rho =-\frac{\sqrt{t(t-2)}}{2t-1}\ (t>2)$.
\end{example}
\begin{example}\label{example5.2}
Let
$$
F(x,y)=1-\frac{1}{(y+1)^t}-\frac{1}{(x+1)^t}+\frac{1}{(x+y+1)^t}\quad (x>0, y>0, t>0)
$$
be a bivariate distribution with marginal distributions  $H(x)=1-\frac{1}{(x+1)^t}\ (x>0)$ and
$G(y)=1-\frac{1}{(y+1)^t}\ (y>0)$.  It follows from (\ref{2.2}) that for any $t>0$
$$
\tau=\frac{1}{2t+1}.
$$ 
A corresponding simulation experiment  made for different values of $t>0$ confirmed this result.

Observe that here $\rho =\frac{1}{t}\ (t>2)$.
\end{example}
\begin{example}\label{example5.3}
Let
$$
F(x,y)=xy(1+\alpha (1-x)(1-y))\quad (0<x,y<1, -1\leq \alpha \leq 1)
$$
be a bivariate distribution with marginal distributions  $H(x)=x\ (0<x<1)$ and
$G(y)=y\ (0<y<1)$.  It follows from (\ref{2.2}) that for any $\alpha $
$$
\tau=\frac{2\alpha }{9}.
$$ 
A corresponding simulation experiment  made for different values of $\alpha $ confirmed this result.

Observe that here $\rho =\frac{\alpha }{3}\ (-1\leq \alpha \leq 1)$.
\end{example}
The paper is submitted to a statistical journal.
\section*{References}
\begin{description} {\small
\item Bairamov, I. and Stepanov, A. (2010).\ Numbers of
near-maxima for the bivariate case,  {\it Statistics  $\&$ Probability
Letters}, {\bf 80}, 196--205.

\item  Bhattacharya, B.B. (1974).\ Convergence of sample paths of
normalized sums of induced order statistics, {\it Ann. Statist. }, {\bf 2}, 1034--1039.

\item Bhattacharya, B.B. (1984).\ Induced order statistics: theory and applications, In {\it Handbook of Statistics 4}, Ed.
Krishnaiah, P. R. andSen, P. K.,  North Holland, Amsterdam, 383--403. 

\item Chu, S.J.,  Huang, W.J. and Chen, H. (1999).\ A study of
asymptotic distributions of concomitants of certain order statistics, {\it Statist. Sinica}, {\bf 9}, 811--830.

\item Daniels, H. E. (1950).\ Rank correlation and population models, {\it Journal of the Royal Statistical Society}, Ser. B, {\bf 12} (2),
171–-191.

\item David, H.A. (1994).\ Concomitants of Extreme Order Statistics, In {\it Extreme Value Theory and Applications},
Proceedings of the Conference on Extreme Value Theory and Applications, {\bf 1}, Ed. Galambos, J., Lechner, J., and Simiu,
E., Kluwer Academic Publishers, Boston  211--224.

\item David, H.A. and Galambos, J. (1974).\ The asymptotic theory of concomitants of order statistics, {\it J. Appl. Probab.}, {\bf 11}, 762--770.

\item David, H.A. and Nagaraja, H.N. (2003).\  {\it Order Statistics}, Third  edition, John Wiley \& Sons, NY.

\item Durbin J. and Stuart,  A. (1951).\ Inversions and rank correlation coefficients, {\it Journal of the Royal Statistical Society},  Ser. B {\bf 13} (2), 303–-309.

\item Egorov, V. A. and Nevzorov, V. B. (1984).\ Rate of
convergence to the Normal law of sums of induced order statistics,
{\it Journal of Soviet Mathematics} (New York), {\bf 25}, 1139--1146.

\item Fisher, R. A. (1921).\ On the 'probable error' of a coefficient of correlation deduced from a small sample, Metron, {\bf 1},  3–-32.

\item Goel, P. K. and Hall, P. (1994).\ On the average difference between concomitants and order statistics,
{\it  Ann. Probab.}, {\bf 22}, 126--144.

\item Kendall, M. G. (1970).\ {\it Rank Correlation Methods}, London, Griffin.

\item   Xu, W., Hou,  Y.,   Hung, Y. S. and Zou, Y. (2009). Comparison of Spearman's rho and Kendall's
tau in Normal and Contaminated Normal Models, {\it arXiv:1011.2009v1 [cs.IT]}.

 }
\end{description}
\end{document}